\magnification=1200
\advance\hoffset by 1,5 truecm
\advance\hsize   by -2,3 truecm
\overfullrule=0mm
\parindent=0cm

%--------- Lettres cursives  

\font\timbfXVIII=cmbx10 at 18pt 
\font\timbfXIV=cmbx10 at 14pt 
\font\timbfXII=cmbx10 at 12pt 

\def\P{\hbox{{\rm P}\kern-.85em\hbox{{\rm I}}}\ }
\def\F{\hbox{{\rm F}\kern-.8em\hbox{{\rm I}}}\ }
\def\K{\hbox{{\rm K}\kern-.9em\hbox{{\rm I}}}\ }
\def\H{\hbox{{\rm H}\kern-.9em\hbox{{\rm I}}}\ }
\def\N{\hbox{{\rm N}\kern-.9em\hbox{{\rm I}}}\ \ }
\def\Z{\hbox{{\rm Z}\kern-.4em\hbox{{\rm Z}}}\ \!}
\def\Q{\hbox{{\rm Q}\kern-.55em\hbox{\vrule width .5pt height 6.5pt }}\ \ }
\def\R{\hbox{{\rm R}\kern-.9em\hbox{{\rm I}}}\ \ }
\def\B{\hbox{{\rm B}\kern-.9em\hbox{{\rm I}}}\ \ } 
\def\C{\hbox{{\rm C}\kern-.45em\hbox{\vrule width .5pt height 6.5pt}}\ \ }
\def\D{\hbox{{\rm D}\kern-.45em\hbox{\vrule width .5pt height 6.5pt}}\ \ }
\def\U{\hbox{{\rm U}\kern-.9em\hbox{{\rm I}}}\ }
\def\A{\hbox{{\rm A}\kern-.88em\hbox{{\it I}}}\ }
%*********************************************************

%Pour avoir un grand ¤ :TAPER \SS 
\font\grandsy=cmsy10 scaled \magstep2
\def\SS{{\grandsy x}}

%Pour avoir un petit ¤ :TAPER \SSpetit 

\centerline{\timbfXVIII STONE-WEIERSTRASS THEOREM}
\bigskip
\centerline{\timbfXIV G. LAVILLE and I.P. RAMADANOFF}
\bigskip
\centerline{\vbox{\hrule width 1cm}}

\vskip 1cm
{\bf Abstract.-} It will be shown that the Stone-Weierstrass theorem for Clifford-valued
functions is true for the case of even dimension. It remains valid for the odd dimension if
we add a stability condition by principal automorphism.

\bigskip\bigskip
{\bf Introduction.-} \ Recall the classical Stone-Weierstrass theorem :  let $Y$ be a metric
space, \ ${\cal C} (Y ; \R )$  the set of all continuous functions from $Y$ in \ $\R$, \break 
$B\subset {\cal C} (Y ; \R )$  a subset such that $B$ contains the constant function 1 and
separates the points of $Y$. Then the algebra  $A_B(Y ; \R )$, generated by $B$ is dense in
${\cal C} (Y ; \R )$  for the topology of the uniform convergence on every compact.
\medskip
\hskip 0,5cm  It is well-known that if one substitutes the field \ $\R$  by $\C$, then an
additional hypothesis is needed, namely :  $B$ should be stable with respect to complex
conjugation.  In case we are omitting this hypothesis and if we take, for example, $Y$ to be
an open subset of $\C$  and $Y = \{ 1,z\}$, then we will get the algebra of holomorphic
functions.

\medskip
\hskip 0,5cm  Let us mention that the case of functions taking values in the quaternonian
field is known [2] and it is analogous to the real case.

\medskip
\hskip 0,5cm  Here, we will investigate the situation when  \ $\R$  is replaced by \
$\R_{p,q}$ - an universal Clifford algebra of \ $\R^n$, \ $n=p+q$,  with a quadratic form of
signature  $(p,q)$.  This study is motivated by the theory of monogenic functions [1].  The
present paper is organized as follows : in the  \SS 1 we will recall some notations usually
employed in Clifford algebras. The \SS 2 will deal with some elements of combinatorics. The
essential part of the paper in the \SS 3 in which we give a formula allowing to compute the
scalar part of a given Clifford number. As an application of this formula, we are able to
prove in \SS 4 the following Stone-Weierstrass theorem for ${\cal C} (Y ; \R_{p,q})$:

\bigskip\bigskip
{\bf Theorem.-} \ {\it  Let $Y$ be a metric space and ${\cal C} (Y ; \R_{p,q})$  the set of
all continuous functions from $Y$ to $\R_{p,q}$.  Let $B\subset {\cal C} (Y, \ \R_{p,q})$  be
such that $B$ contains the constant function 1 and separates the points of  $Y$. If $p+q$ is
odd, suppose in addition that $B$ is stable with respect to the principal automorphism $\ast$
.  Then, the algebra  $A_B (Y ; \ \R_{p,q})$, generated by $B$, is dense in ${\cal C} (Y ; \
\R_{p,q})$  for the topology of uniform convergence on compact sets.}

\bigskip\bigskip
{\timbfXII 1. Notations}
\medskip
\hskip 0,5cm In a Clifford algebra \ $\R_{p,q} = C_0 \oplus C_1 \oplus \ldots \oplus
C_n$,with $n=p+q$,  the spaces $C_0, C_1, \ldots , C_n$  are supposed to be of respective
basis  $\{ 1\}$, $\{ e_1, e_2,\ldots ,e_n\}$,

 ${\{ e_{ij}\}_{i<j}, \ldots , \{ e_{i_{1} \ldots i_{k}}
\}}_{i_{1}<i_{2}<\cdots < i_{k}}, \ldots , \{ e_{1.2\ldots n}\}$, where $(i_1,\ldots , i_k)$ 
is a multiindex with $i_1,\ldots , i_k\in \{ 1,\ldots , n\}$, $1\leq i_1 < \ldots < i_k \leq
n$.  The algebra obeys to the laws :

$$\cases{ 
e_i^2 = 1,                                                  &$i=1,\ldots , p$\cr
e_i^2 = -1,                                                 &$i=p+1,\ldots ,n$\cr
e_ie_j = -e_je_i,                                           &$i\not= j$\cr
e_{i_{1}\ldots i_{k}} = e_{i_{1}}e_{i_{2}}\cdots e_{i_{k}}, &for $i_1<i_2<\ldots < i_k$\cr}$$

We will make use of the decomposition of a Clifford number a in its scalar (real) part
$\langle a\rangle_0$, its 1-vector $\langle a\rangle_1 \in C_1$, its bivector part $\langle
a\rangle_2 \in C_2$, etc $\ldots$ up to its pseudo-scalar part $\langle a\rangle_n \in C_n$,
i.e :
$$a = \langle a \rangle_0 + \langle a\rangle_1 +\cdots + \langle a\rangle_n,$$
where : 
$$\langle a\rangle_k = \sum_{\scriptstyle J\atop\scriptstyle 
\mid J\mid  = k}\  a_J \ e_J.$$

Where $J = (j_1,\ldots , j_k)$  is a multiindice and $\mid J\mid = k$, \ \ $e_J = e_{j_{1}}
\cdots e_{j_{k}}$.

Recall that the principal involution $_\ast$ , the anti-involution $^\ast$  and the reversion
$\sim$ act on $a\in \R_{0,n}$  as follows :
$$\eqalign{
&a_\ast = \sum_{k=0}^n \  (-1)^k \langle a \rangle_k \cr
&a^\ast = \sum_{k=0}^n \ (-1)^{k(k+1)\over 2} \langle a\rangle_k \cr
&a^\sim = \sum_{k=0}^n \ (-1)^{k(k-1)\over 2} \ \langle a\rangle_k \cr}$$

Now, define
$$e^i  = \cases{
e_i,  &if \ \ $1\leq i\leq p$\cr
-e_i, &if \ \ $p+1 \leq i\leq p+q$\cr}$$

and $e^J = e^{j_{k}} \cdots e^{j_{1}}$.

\eject
{\timbfXII2. Some combinatorics}
\medskip
\hskip 0,5cm Let us study the partition of the set $\{ 1,\ldots , n\}$  in two strictly ordered
subsets: $I = \{ i_1,\ldots ,i_k\}$ and $J = \{ j_1,\ldots ,j_p\}$.  As for as the relative
position of $J$ with respect to $I$ is concerned, we have different possible cases : $J\cap I
= \phi$~; just one $j_\alpha$ belongs to $I ; \ldots ; \ell$ among the $j_\alpha 's$ belong to
$I ; \ldots ;$  the largest  possible number of $j_\alpha 's$ belongs to $I$. It is easy to
compute the cardinals of the corresponding sets : 

For the first case, the cardinal is  $C_{n-k}^p \ C_k^{\sup\{ 0, p-(n-k)\} }$.  If just one
$j_\alpha$ belongs to $I$, then we will have  $C_{n-k}^{p-1} \ C_k^{\sup \{ 0, p-(n-k)\} +1}$ 
and so on $\ldots$ \ In the last case, we will get  $C_{n-k}^0 \ C_k^{\inf \{ p,k\} }$.

\medskip
\hskip 0,5cm Now, recall the following result which is well-known in classical probability
theorey [3] :

\bigskip\bigskip
{\bf Lemma 1.-} \ {\it  For every  $k, \ 0\leq k \leq n$ :
$$\sum_{\ell = \sup \{ 0, p-(n-k)\} }^{\inf \{ p,k\} } \ C_{n-k}^{p-\ell} \ C_k^\ell = C_n^p
.$$

In fact, this lemma will not be used here, but its elementary proof, which will be given
below, is a source of inspiration for the next result  (Lemma 2).
}

\bigskip
{\bf Proof of Lemma 1 --} For every $k$, \ $0\leq k\leq n$, one has \break 
$(1+x)^{n-k} (1+x)^k = (1+x)^n$,  which involves  

$\displaystyle\sum_{\ell = 0}^k \ (1+x)^{n-k}\ C_k^\ell \ x^\ell = \displaystyle\sum_{p=0}^n \
C_n^p \ x^p$, and again : 

$\displaystyle\sum_{\ell = 0}^k \ \displaystyle\sum_{n=0}^{n-k} \
C_{n-k}^n \ x^n \ \ C_k^\ell \ x^\ell = \displaystyle\sum_{p=0}^n \ C_n^p x^p$.  Let us set
$n+\ell = p$, i.e. $n=p-\ell$. 
\smallskip
Then the double sum is equal to 

$\displaystyle\sum_{\ell =
0}^k \ \ \displaystyle\sum_{p=\ell}^{n-k+\ell} \ C_{n-k}^{p-\ell} \ C_k^\ell x^p =
 \displaystyle\sum_{p=0}^n \ \ \displaystyle\sum_{\ell = \sup \{ 0,p -(n-k)\} }^{\inf \{ p,k\}
} \ C_{n-k}^{p-\ell} \ C_k^\ell \  x^p$.  \hfill$\diamondsuit$

\bigskip
\hskip 0,5cm It just remains to indentify the coefficients of  $x^p$. Now, we are in a
position to formulate and prove the following :

\bigskip
{\bf Lemma 2.-} \ {\it 
$$\sum_{p=0}^n \ \sum_{\ell = \sup \{ 0,p-(n-k)\} }^{\inf \{ p,k\} } \ (-1)^{pk+\ell} \
C_{n-k}^{p-\ell} \ C_k^\ell = \cases{
0, &if \ \ $1\leq k\leq n-1$\cr
0, &if \ \ $k=n, n$ even \cr
2^n, &if \ \ $k=n,n$ odd\cr
2^n, &if \ \ $k=0$.\cr}$$
}

\bigskip
{\bf Proof of Lemma 2 --}~ Start from
$(1+(-1)^kx)^{n-k} (1+(-1)^{k+1}x)^k = $
$$\eqalign{ &= \sum_{\ell = 0}^k (1+(-1)^kx)^{n-k} (-1)^{(k+1)\ell } C_k^\ell x^\ell = \cr
&=\sum_{\ell = 0}^k \sum_{n=0}^{n-k} (-1)^{kn}\  C_{n-k}^n\  x^n (-1)^{(k+1)\ell}\  C_k^\ell
\ x^\ell = \cr
&=\sum_{p=0}^n \ \sum_{\ell = \sup \{ 0,p-(n-k)\} }^{\inf \{ p,k\} } (-1)^{pk+\ell}
\ C_{n-k}^{p-\ell} \ \ C_k^\ell \ x^p , \cr}$$

because  $kn+(k+1)\ell = pk+\ell$. Thus it is enough to set $x=1$ and remark that :
$$(1+(-1)^k)^{n-k} (1+(-1)^{k+1})^k = \cases{
2^n, &if \ \ $k=0$\cr
0,   &if \ \ $1\leq k \leq n-1$\cr
2^n, &if \ \ $k=n, n$ \ odd\cr
0,   &if \ \ $k=n, \ n$ \ even\cr}$$ \hfill$\diamondsuit$

\vskip 1cm
{\timbfXII{3. A formula for the real part of $a\in \R_{p,q}$ }}
\medskip
{\bf Lemma 3.-} \ {\it For every multiindice $J$, we have $e_J \ e^J = 1$.
}

\bigskip
{\bf Lemma 4.-} \ {\it  Let $I = (i_1,\ldots ,i_k)$, \ $\mid I\mid = k$. \ $J = (j_1,\ldots ,
j_p)$, \ $\mid J\mid  = p$  there is the following equality
$$\sum_{p=0}^n \ \sum_{\mid J \mid =p} \ e_Je_Ie^J = \cases
{2^n &if $k=0$ or if $k=n$ with $n$ odd \cr
0 &in other cases \cr}$$
}

\bigskip
{\bf Proof} --  Decompose the sum 
$$\sum_{\mid J\mid =p} \ e_je_I e^J$$
following the relative position of $J$ with respect to $I$. If  $J\cap I = \phi$  we have 
$C_{n-k}^p C_k^0$  such possibilities and the anticommutation gives  $(-1)^{pk}$.  

If only one $j_\alpha\in I$ we have  $C_{n-k}^{p-1} \ C_k^1$  such possibilites and the
anticommutation gives  $(-1)^{(p-1)k} \ (-1)^{k-1}$ and so on, $\ldots$ , if \ $\ell \
j_\alpha\in I$  we have  $C_{n-k}^{(p-\ell )k} C_k^\ell$  such possibilities and the
commutation gives  $(-1)^{(p-\ell )k} \ (-1)^{\ell (k-1)}$.

The sum is equal to
$$\sum_{\ell = \sup \{ 0, p-(n-k)\} }^{\inf \{ p,k)} \ (-1)^{(p-\ell )k} \ (-1)^{\ell (k-1)}
\ C_{n-k}^{p-\ell} \ C_k^\ell e_I$$

Thus we could apply lemma 2 and the result follows.  \hfill$\diamondsuit$

\bigskip
\hskip 0,5cm  The next result is a formula for the scalar part of a Clifford number.

\bigskip\bigskip
{\bf Theorem 1.-} \ {\it  Let $a\in \R_{p,q}$. Then : 

{\parindent=0,6cm
\item{\bf a)} if $n$ is even,
$$\langle a\rangle_0 = {1\over 2^n} \ \sum_{p=0}^n \ \sum_{\mid J\mid =p} e_J ae^J.$$

\item{\bf b)} if $n$  is odd,
$$\langle a \rangle_0 = {1\over 2^{n+1}} \ \sum_{p=0}^n \ \sum_{\mid J\mid =p} \ e_Ja \ e^J +
{1\over 2^{n+1}} \ \sum_{p=0}^n \ \sum_{\mid J\mid =p} \ e_J a_\ast\  e^J.$$
\par}
}

\bigskip
{\bf Proof} --  When $a\in \R_{0,n}$,  then  
$$a = \sum_{k=0}^n \ \sum_{\mid I\mid =k} a_I e_I,$$
where  $I = (i_1,\ldots ,i_k)$, \ $1\leq i_1 < i_2 < \ldots < i_k \leq n$.  Take the sum 
$$\sum_{p=0}^n \ \sum_{\mid J\mid =p} \ e_J a \ e^J = \sum_J \ \sum_I \ a_I\  e_J \ e_I\  e^J.$$

Now, apply lemma 4 :

{\parindent=0,6cm
\item{\bf a)} if  $n$  is even, one gets :
$$\sum_{p=0}^n \ \sum_{\mid J\mid =p} \ e_J a\  e^J = 2^n \ \langle a \rangle_0,$$

\item{\bf b)} if $n$ is odd, one has :
$$\sum_{p=0}^n \ \sum_{\mid J\mid =p} \ e_J a \ e^J = 2^n \ \langle a\rangle_0 + 2^n \ \langle a
\rangle_n .$$
\par}

But, in the case when $n$ is odd, \ $\langle a_\ast\rangle_n = (-1)^n \ \langle a\rangle_n =
- \langle a\rangle_n$.  Thus, we get the part b)  of the theorem. \hfill$\diamondsuit$

\bigskip
{\bf Remark.-}  For $n=1$, the preceding formula becomes to  

$4Re \ a = (a-iai)+(\overline a - i\overline ai)$ in $\R_{0,1} = \C$  with the classical
notations of \ $\C$.

\medskip
\hskip 0,5cm  For $n=2$,  this means that  $4Re \ a = a - iai - jaj - kak$ in \ $\R_{0,2} =
\H$  with the classical notations of \ $\H$, \ [2].

\vskip 1cm
{\timbfXII{4. The Stone-Weierstrass theorem for ${\cal C} (Y ; \ \R_{p,q})$.}}
\medskip
{\bf Theorem 3.-} \ {\it  Let $Y$ be a metric space and  ${\cal C}\ (Y ; \ \R_{p,q})$ the set
of continuous functions from $Y$ into \ $\R_{p,q}$. Let $B\subset {\cal C} (Y ; \ \R_{p,q})$ 
be such that $B$  contains the constant function 1 and separates the points of  $Y$.  When
$p+q$ is even, nothing more is supposed. If $p+q$ is odd, suppose $B$ be stable with respect
to the principal involution $\ast$.

\medskip
\hskip 0,5cm  Then, the algebra  $A_B(Y ; \ \R_{p,q})$, generated by $B$, is dense in ${\cal
C} (Y ; \ \R_{p,q})$  for the topology of uniform convergence on compact.}

\bigskip
{\bf Proof} --  Set $A_B(Y ; \ \R )$ for the subspace of  $A_B(Y ; \ \R_{p,q})$  consis\-ting of those
functions which take real values. This is a real algebra. Let  $A_B(Y ; \ \R )_I$  be the
subspace of  $A_B(Y ; \ \R_{p,q})$  consisting of the $I$-components of functions from 
$A_B(Y ; \ \R_{p,q})$.  Thus, we have   $f_I = \langle f \ e^I\rangle_0$ and $A_B(Y ; \ \R
)_I \subset A_B (Y ; \ \R )$  by theorem 2.

In this way,  $A_B(Y ; \ \R )$  satisfies to the hypothesis of the classical Stone-Weierstrass
theorem for real functions. The algebra  $A_B(Y ; \ \R )$  is consequently dense in ${\cal C}
(Y ; \ \R )$. Finally, one can conclude that :
$$A_B(Y ; \ \R_{p,q}) = \bigoplus_I \ A_B(Y ; \ \R ) e_I$$
is dense in ${\cal C} (Y ; \ \R_{p,q})$. \hfill$\diamondsuit$

\vskip 1cm
{\timbfXII{5. A remark}}
\medskip
\hskip 0,5cm It should be noted that the computations of the scalar part it strongly related
to formulas related to the Hestenes multivector derivative : see [4], chapter 2.

%%%%%%%%%%%%%%%%%%%%%%%%%%%%%%%%%
% Pour taper une bibliographie
%********************************************
\newdimen\margeg \margeg=0pt
\def\bb#1&#2&#3&#4&#5&{\par{\parindent=0pt
    \advance\margeg by 1.1truecm\leftskip=\margeg
    {\everypar{\leftskip=\margeg}\smallbreak\noindent
    \hbox to 0pt{\hss\bf [#1]~~}{\bf #2 - }#3~; {\it #4.}\par\medskip
    #5 }
\medskip}}

%Autre Modele avec le nom de l'auteur en sc
\newdimen\margeg \margeg=0pt
\def\bbaa#1&#2&#3&#4&#5&{\par{\parindent=0pt
    \advance\margeg by 1.1truecm\leftskip=\margeg
    {\everypar{\leftskip=\margeg}\smallbreak\noindent
    \hbox to 0pt{\hss [#1]~~}{\pmb{\sc #2} - }#3~; {\it #4.}\par\medskip
    #5 }
\medskip}}

\eject
{\bf References :}
\bigskip

\bb 1&R. DELANGHE, F. SOMMEN, V. SANC\v EK&Clifford Algebra and Spinor-valued functions&
         Kluwer& &

\bb 2&J. DUGUNDJI&Topology&Allyn and Bacon & &

 \bb 3&W. FELLER&An introduction to the theory of Probability and its applications&J. Wiley& &

\bb 4&D. HESTENES, G. SOBCZYK&Clifford Algebra to Geometric Calculus &Reidel  & &

\bigskip\bigskip
\hfill{\bf Universit\'e de Caen}

\hfill{\sl  D\'epartement de Math\'ematiques}

\hfill{\sl  Esplanade de la Paix}

\hfill {\sl 14032 CAEN CEDEX - FRANCE}

\end